\documentclass[preprint,11pt]{article}
\usepackage[top=2.5cm,bottom=2.5cm,left=2.5cm,right=2.5cm]{geometry}

\usepackage{cite}
\usepackage{enumitem}
\usepackage{epsfig}
\usepackage{graphicx}
\usepackage{caption}
\usepackage{multirow}
\usepackage{dcolumn}
\usepackage{bm}
\usepackage{amsmath,amsfonts,amssymb,amsthm,amscd,bbm,mathrsfs}
\usepackage{colortbl}
\usepackage{booktabs}
\usepackage{alltt}
\usepackage{moreverb}
\usepackage[]{algorithm,algcompatible,lipsum}
\usepackage[]{algpseudocode}
\usepackage{listings}
\usepackage[dvipsnames,table]{xcolor}
\usepackage{xurl}
\usepackage{hyperref}
\hypersetup{colorlinks,citecolor=black,filecolor=black,linkcolor=black,
	urlcolor=black,bookmarksnumbered}

\definecolor{lightblue}{RGB}{230, 240, 255}

\newcommand{\dparc}[4]{\frac{\partial #1}{\partial #2 \partial #3
		\partial #4}}

\newcommand{\dparn}[5]{\frac{\partial #1}{\partial #2 \partial #3
		#4 \partial #5}}

\definecolor{darkgreen}{rgb}{0.0, 0.5, 0.0} 
\definecolor{darkred}{rgb}{0.5, 0.0, 0.0}

\lstdefinestyle{fortran-style}{
  language=[90]Fortran,
  basicstyle=\ttfamily\small,
  morekeywords={abstract, import, asinh, acosh , atanh},
  keywordstyle=\color{blue}\bfseries,
  commentstyle=\color{darkgreen},
  stringstyle=\color{darkred},
  numbers=left,
  numberstyle=\tiny\color{gray},
  stepnumber=1,
  numbersep=5pt,
  showspaces=false,
  showstringspaces=false,
  showtabs=false,
  tabsize=4,
  captionpos=b,
  breaklines=true,
  breakatwhitespace=false,
  escapeinside={\%*}{*)}
}

\lstdefinelanguage{Mathematica}{
  morekeywords={
    Module,Function,With,Do,Map,Array,Range,Pick,Transpose,Normal,
    ClearAll,StringRiffle,ToString,CoefficientRules,BellY,Symbol
  },
  sensitive=true,
  morecomment=[l]{(*},
  morecomment=[s]{(*}{*)},
  morestring=[b]{"}
}

\makeatletter
\newenvironment{breakablealgorithm}
  {
  \begin{center}
    \refstepcounter{algorithm}
    \hrule height.8pt depth0pt \kern2pt
    \renewcommand{\caption}[2][\relax]{
      {\raggedright\textbf{\fname@algorithm~\thealgorithm} ##2\par}%
      \ifx\relax##1\relax 
      \addcontentsline{loa}{algorithm}{\protect\numberline{\thealgorithm}##2}%
      \else 
      \addcontentsline{loa}{algorithm}{\protect\numberline{\thealgorithm}##1}%
      \fi
      \kern2pt\hrule\kern2pt
    }
}
  {
    \kern2pt\hrule\relax
  \end{center}
  }
\makeatother

  \title{Efficient Computation of Arbitrary-Order Directional
    Derivatives in Multiple Directions via Generalized Dual Numbers}

  \author{
  F. Pe\~nu\~nuri,~
  K. B. Cant\'un-Avila,~
  O. Carvente,~
  R. Pe\'on-Escalante
}

\begin{document}
  \maketitle
  \vspace*{-0.7cm}
\begin{center}
\footnotesize{Facultad de Ingenier\'ia, Universidad
Aut\'onoma de Yucat\'an, A.P. 150, Cordemex, M\'erida, Yucat\'an,
M\'exico.}
\end{center}

  \begin{abstract}
    Arbitrary-order directional derivatives along multiple (possibly
    distinct) directions are computed through a generalized dual-number
    formulation for both scalar- and vector-valued functions. The proposed
    framework combines generalized dual evaluations with an
    inclusion--exclusion reconstruction of symmetric multilinear forms,
    allowing general multidirectional derivatives to be reconstructed from
    repeated-direction evaluations without explicitly constructing
    higher-order derivative tensors. Mixed directional and mixed partial
    derivatives arise naturally as particular cases of the formulation.
    The methodology further enables the systematic computation of
    arbitrary-order kinematic quantities and the construction of
    Taylor-series methods for systems of ordinary differential equations
    through automatically generated time derivatives. Numerical examples
    include the computation of high-order directional derivatives in
    high-dimensional functions, mixed partial derivatives, arbitrary-order
    kinematic quantities, and Taylor-series integration methods. The
    implementation is developed in modern \textsc{Fortran} within an
    open-source framework compatible with the Fortran Package Manager
    ecosystem.
  \end{abstract}

  \textbf{Keywords:} Directional derivatives, dual numbers, automatic
  differentiation, kinematic quantities, Taylor series methods.
  
\vspace*{0.7cm}

\section{Introduction}\label{sec1}
Directional derivatives of arbitrary order along multiple distinct
directions arise naturally in multilinear analysis, high-order Taylor
methods, sensitivity analysis, uncertainty propagation, and the
kinematic modeling of mechanical systems
\cite{GriewankWalther2008,MartinsHwang2013,HairerLubichWanner2006,
LeimkuhlerReich2004,Cantun2025}. While first- and second-order
derivatives are sufficient for many applications, several emerging
problems involving higher-order effects require more general
differentiation frameworks capable of systematically handling
arbitrary-order derivatives.

Traditional approaches for computing such derivatives become rapidly
impractical as the order increases. Symbolic differentiation is often
computationally expensive and difficult to manage for large expressions
\cite{Atilim2018,Peon2020}, whereas finite-difference approximations are
affected by truncation and cancellation errors
\cite{Squire1998,Martins2003}. Automatic differentiation (AD) provides
an efficient alternative by generating derivatives directly from the
computational structure of a function to machine precision
\cite{Griewank1989,Bischof1992,CHINCHALKAR1994197,GAWA2003}. In
particular, formulations based on dual numbers have proved especially
convenient because derivatives can be obtained through algebraic
extensions of real arithmetic
\cite{Cheng1994,Jefrey2011,Kalos2021,Peon2025}.

Existing generalized dual-number approaches efficiently compute
directional derivatives along a repeated direction,
$d_n(\mathbf{v},\mathbf{v},\ldots,\mathbf{v})$,
through a single generalized dual evaluation
\cite{PEONESCALANTE2024,Penunuri2025}. However, many applications
require the more general multidirectional derivative
$d_n(\mathbf{v}_1,\mathbf{v}_2,\ldots,\mathbf{v}_n)$,
where the directions are not necessarily equal. Such derivatives arise
naturally in multilinear expansions, mixed sensitivity analysis,
tensor reconstruction, and mixed partial differentiation. Although they
can in principle be recovered from higher-order derivative tensors,
explicit tensor construction rapidly becomes computationally prohibitive
as the order of differentiation increases.

The present work addresses this problem by combining generalized
dual-number evaluations with an inclusion--exclusion reconstruction of
symmetric multilinear forms. This permits the computation of
arbitrary-order directional derivatives along multiple distinct
directions without explicitly constructing higher-order derivative
tensors. Mixed partial derivatives arise as a particular case through
the use of canonical basis vectors.

The proposed framework applies to both scalar-valued functions
$f:\mathbb{R}^m\to\mathbb{R}$ and vector-valued functions
$\mathbf{f}:\mathbb{R}^m\to\mathbb{R}^p$. Beyond the recovery of mixed
partial derivatives, the methodology enables the systematic computation
of kinematic quantities of arbitrary order, extending beyond the
velocity, acceleration, jerk, and snap quantities commonly considered
in the literature. In addition, the proposed formulation naturally
leads to Taylor-series methods for systems of ordinary differential
equations through the automatic generation of higher-order time
derivatives.

The computational implementation is developed in modern
\textsc{Fortran} within the DNAOAD
infrastructure~\cite{Penunuri2025}, with modifications aimed at
facilitating compilation and execution through the Fortran Package
Manager ecosystem. The implementation used in this work has been
archived with a persistent identifier on
Zenodo~\cite{penunuri_2026_zDNAODD}, ensuring reproducibility of the
reported results.

\section{Dual numbers and arbitrary-order directional derivatives}

The algebraic framework used in this work builds on the dual-number
formulation for higher-order directional derivatives developed
in~\cite{PEONESCALANTE2024} and on its arbitrary-order implementation
reported in~\cite{Penunuri2025}. In that setting, a generalized dual
number is written as
\begin{align}
  r = \sum_{k=0}^n r_k \epsilon_k,
\end{align}
where the basis elements satisfy
\begin{equation}\label{gentabmult}
  \epsilon_i \epsilon_j =
  \begin{cases}
    0, & i+j>n,\\[0.2cm]
    \dfrac{(i+j)!}{i!j!}\epsilon_{i+j}, & \text{otherwise}.
  \end{cases}
\end{equation}

This algebra permits the evaluation of higher-order directional
derivatives by replacing the real argument of a function with a
generalized dual argument. In particular, for an analytic function
$f:\mathbb{R}^m\to\mathbb{R}$, a point
$\mathbf{q}\in\mathbb{R}^m$, and a direction
$\mathbf{v}\in\mathbb{R}^m$, the evaluation
$f(\mathbf{q}+\epsilon_1\mathbf{v})$ contains, in its
$\epsilon_k$ component, the $k$th-order directional derivative along
the repeated direction $\mathbf{v}$; that is,
\begin{align}\label{repeated_direction_derivative}
  f(\mathbf{q}+\epsilon_1\mathbf{v}).\epsilon_k
  =
  d_{kf\mathbf{q}}(\mathbf{v},\ldots,\mathbf{v}),
  \qquad k=1,\ldots,n.
\end{align}
Here $f.\epsilon_k$ denotes the coefficient of $\epsilon_k$ in the
generalized dual expansion. The derivation and computational details of
this construction are given in~\cite{PEONESCALANTE2024,Penunuri2025}.

The previous formulation directly provides derivatives along a single
repeated direction. The objective here is different: we seek the general
symmetric multilinear form
\begin{align}
  d_{nf\mathbf{q}}(\mathbf{x}_1,\mathbf{x}_2,\ldots,\mathbf{x}_n),
\end{align}
where the directions $\mathbf{x}_1,\ldots,\mathbf{x}_n$ are not
necessarily equal. This generalization is essential for computing mixed
directional derivatives and, as a particular case, mixed partial
derivatives of arbitrary order.

For compactness, define
\begin{align}
  d_{nf\mathbf{q}}(\mathbf{v})
  :=
  d_{nf\mathbf{q}}(\mathbf{v},\mathbf{v},\ldots,\mathbf{v}).
\end{align}
Since $d_{nf\mathbf{q}}$ is a symmetric multilinear form, evaluating it
at the sum of $n$ directions gives
\begin{equation}\label{sum_multilinear_expansion}
  d_{nf\mathbf{q}}(\mathbf{x}_1+\mathbf{x}_2+\cdots+\mathbf{x}_n)
  =
  \sum_{i_1,\ldots,i_n\in\{1,\ldots,n\}}
  d_{nf\mathbf{q}}
  (\mathbf{x}_{i_1},\ldots,\mathbf{x}_{i_n}).
\end{equation}
The desired term
$d_{nf\mathbf{q}}(\mathbf{x}_1,\ldots,\mathbf{x}_n)$ can therefore be
isolated by inclusion--exclusion~\cite{BeelerChapter7}. This yields
\begin{align}\label{dnx123n}
  d_{nf\mathbf{q}}(\mathbf{x}_1,\mathbf{x}_2,\ldots,\mathbf{x}_n)
  =
  \frac{1}{n!}
  \sum_{k=1}^{n}(-1)^{n-k}
  \sum_{s_j\in S_k}
  d_{nf\mathbf{q}}\!\left(
    \sum_{\mathbf{x}_p\in s_j}\mathbf{x}_p
  \right),
\end{align}
where $S_k$ denotes the set of all subsets of size $k$ of
$\{\mathbf{x}_1,\mathbf{x}_2,\ldots,\mathbf{x}_n\}$.

For example, when $n=2$,
\begin{align}
  d_{2f\mathbf{q}}(\mathbf{x}_1,\mathbf{x}_2)
  =
  \frac{1}{2}
  \left[
  d_{2f\mathbf{q}}(\mathbf{x}_1+\mathbf{x}_2)
  -d_{2f\mathbf{q}}(\mathbf{x}_1)
  -d_{2f\mathbf{q}}(\mathbf{x}_2)
  \right].
\end{align}

For $n=3$, writing
$\mathbf{x}_{ij}=\mathbf{x}_i+\mathbf{x}_j$ and
$\mathbf{x}_{123}=\mathbf{x}_1+\mathbf{x}_2+\mathbf{x}_3$, one obtains
\begin{align}
  d_{3f\mathbf{q}}(\mathbf{x}_1,\mathbf{x}_2,\mathbf{x}_3)
  =
  \frac{1}{6}
  \Big[
  &d_{3f\mathbf{q}}(\mathbf{x}_{123})
  -d_{3f\mathbf{q}}(\mathbf{x}_{12})
  -d_{3f\mathbf{q}}(\mathbf{x}_{13})
  -d_{3f\mathbf{q}}(\mathbf{x}_{23}) \notag\\
  &+d_{3f\mathbf{q}}(\mathbf{x}_1)
  +d_{3f\mathbf{q}}(\mathbf{x}_2)
  +d_{3f\mathbf{q}}(\mathbf{x}_3)
  \Big].
\end{align}

Equation~\eqref{dnx123n} is the central construction of the present
work. It reduces the computation of a general $n$th-order directional
derivative along multiple distinct directions to a finite number of
repeated-direction evaluations, each of which can be obtained directly
from Eq.~\eqref{repeated_direction_derivative} using generalized dual
numbers.

Mixed partial derivatives are recovered as a special case by choosing
the directions as canonical basis vectors. Thus,
\begin{align}\label{partial_from_directional}
  \dparn{^n f(\mathbf{q})}{q_{i_1}}{q_{i_2}}{\cdots}{q_{i_n}}
  =
  d_{nf\mathbf{q}}
  (\mathbf{e}_{i_1},\mathbf{e}_{i_2},\ldots,\mathbf{e}_{i_n}),
  \qquad
  i_1,\ldots,i_n\in\{1,\ldots,m\}.
\end{align}
For instance,
\begin{align}\label{ejemplod3}
  \dparc{^3 f(\mathbf{q})}{q_i}{q_j}{q_k}
  =
  d_{3f\mathbf{q}}(\mathbf{e}_i,\mathbf{e}_j,\mathbf{e}_k),
  \qquad
  i,j,k\in\{1,2,\ldots,m\}.
\end{align}

The vector-valued case
$\mathbf{f}:\mathbb{R}^m\to\mathbb{R}^p$ follows componentwise, since
each component of $\mathbf{f}$ defines a scalar-valued symmetric
multilinear form. Therefore,
$\mathbf{d}_{n\mathbf{f}\mathbf{q}}
(\mathbf{x}_1,\ldots,\mathbf{x}_n)$ is obtained by applying
Eq.~\eqref{dnx123n} to each component.

Although Eq.~\eqref{partial_from_directional} provides a direct route
to arbitrary-order mixed partial derivatives, the computation of all
partial derivatives becomes combinatorially expensive for large $n$.
In many applications, however, one does not require the complete
higher-order derivative tensor. Instead, selected directional
derivatives are needed along prescribed directions. In such cases, the
combination of generalized dual numbers with
Eq.~\eqref{dnx123n} provides a practical and efficient alternative to
symbolic differentiation or explicit tensor construction.

\section{Applications}
This section presents several representative applications of the proposed
methodology. All functions and illustrative examples discussed herein are
implemented in the accompanying archived
codebase~\cite{penunuri_2026_zDNAODD}.

\subsection{Computation of directional derivatives}
\label{secSPF}

\begin{figure}[htb]
  \centering
  \includegraphics[scale=0.45]{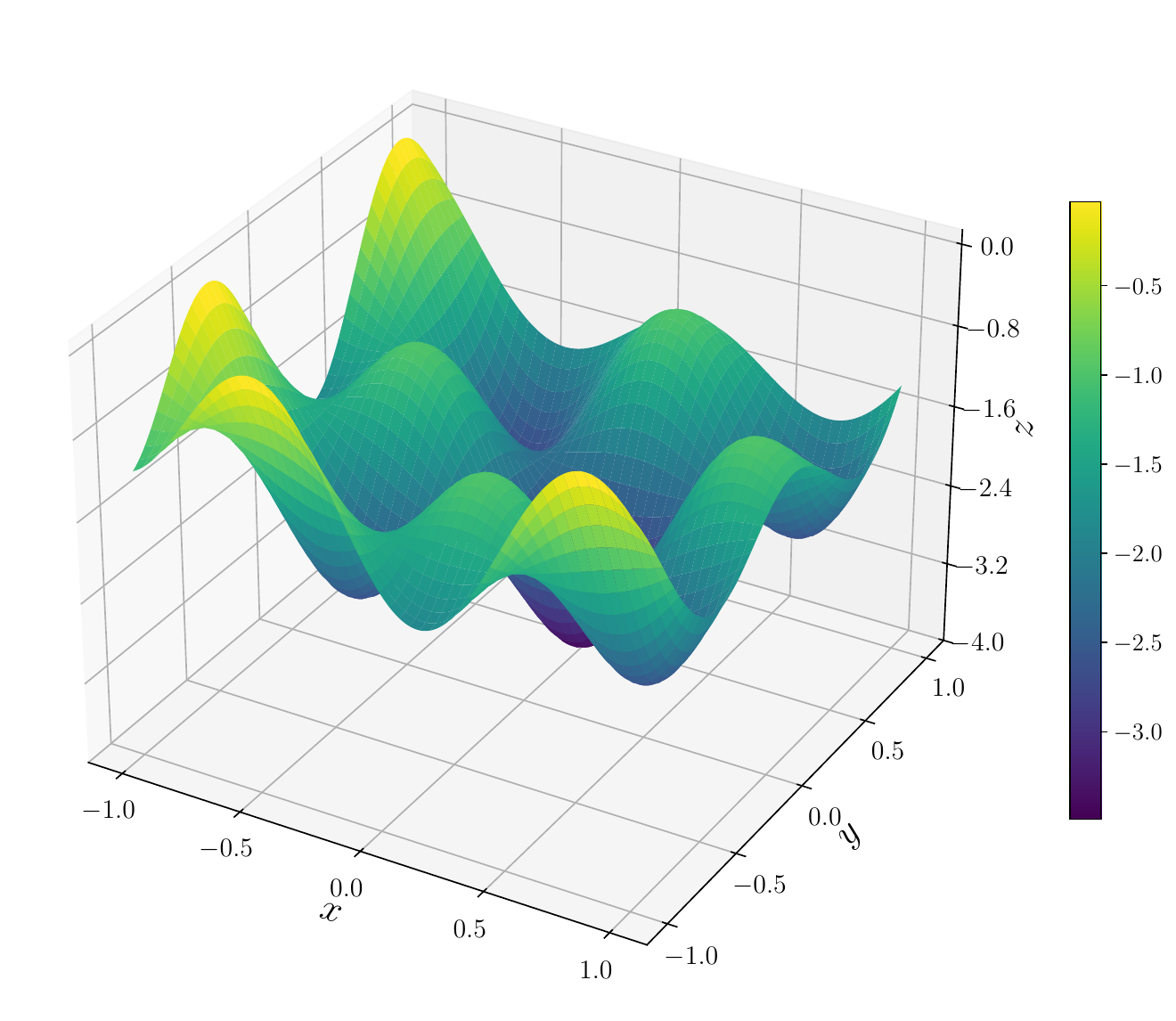}
  \caption{Sinusoidal problem function for two variables.}
  \label{CMFP}
\end{figure}

A classical test function for optimization algorithms is the sinusoidal
problem function~\cite{Penunuri2016}. The $D$-dimensional case is defined as
\begin{align*}
  f_8(\mathbf{q}) &= -A \prod_{i=1}^D \sin(q_i - z)
  - \prod_{i=1}^D \sin[B(q_i - z)], \\
  A &= 2.5, \quad B = 5, \quad z = 30.
\end{align*}
Figure~\ref{CMFP} shows the plot for the two-dimensional case.

We are interested in computing the seventh-order directional derivative
of $f(\mathbf{q})$ along the vector
$\mathbf{x} = [\sin 1, \sin 2, \dots, \sin D]$, evaluated at the point
$\mathbf{q} = [1, 1/2, \dots, 1/D]$, for
$D \in \{100, 1000, 2000, 3000\}$.

Explicitly, the desired derivative is given by
\begin{align}
  d_{7f\mathbf{q}}(\mathbf{x}) &=
  \frac{\partial^7 f(\mathbf{q})}
  {\partial q_{i_1}\partial q_{i_2}\cdots\partial q_{i_7}}
  x_{i_1}x_{i_2}\cdots x_{i_7},
\end{align}

which would be computationally prohibitive if all partial derivatives and
index contractions were evaluated explicitly. In contrast, the proposed
dual-number approach enables efficient computation, with all reported
cases completed in less than 0.1~s on the hardware specified below.
Table~\ref{derdirtab1} reports the computed directional derivatives
together with the corresponding elapsed times.

All computations were performed on a 12th-generation Intel(R) Core(TM)
i7--12650H CPU @ 4.70~GHz with 15~GB of RAM, running GNU/Linux
(AlmaLinux~10) and using the \texttt{gfortran} compiler.

\begin{table}[hbt]
  \caption{Directional derivatives and elapsed times for the sinusoidal
  problem function.}
  \centering
  \scalebox{0.9}{
    \begin{tabular}{rrr}
      \toprule
      $m$ & \multicolumn{1}{c}{Directional derivative} &
      \multicolumn{1}{c}{Time (s)} \\
      \midrule
      \rowcolors{2}{lightblue}{white}
      100  & -12294759.73110 & 0.003 \\
      1000 & -328775.11848   & 0.031 \\
      2000 & -15.34040       & 0.061 \\
      3000 & -0.00032        & 0.092 \\
      \bottomrule
    \end{tabular}
  }
  \label{derdirtab1}
\end{table}

\subsection{Computation of partial derivatives}
The proposed formulation is not primarily intended for the computation
of partial derivatives, since beyond the calculation of Hessians such an
approach may become inefficient. By contrast, quantities of the form
\begin{align}
  \langle \mathbf{u} \mid \mathbf{H}f(\mathbf{q}) \mid \mathbf{v} \rangle
  = \mathbf{u}^{\text{T}} \mathbf{H}f(\mathbf{q}) \mathbf{v},
\end{align}
which represent second-order directional derivatives of
$f(\mathbf{q})$ along the vectors $\mathbf{u}$ and $\mathbf{v}$, can be
computed efficiently. Partial derivatives can nevertheless be recovered
by evaluating directional derivatives along an appropriate vector basis.
The computation of the full Hessian matrix is illustrated in the
algorithm below.

\begin{breakablealgorithm}
  \caption{Pseudocode for \textsc{Hessian}$(f,\mathbf{q})$}
  \label{alg:hessian}
  \begin{algorithmic}[1]
    \Require Procedure $f$ (scalar function of a dual vector variable),
    $\mathbf{q} \in \mathbb{C}^m$
    \Ensure $\mathbf{H} \in \mathbb{C}^{m \times m}$
    \State $m \gets \text{size}(\mathbf{q})$
    \State $\mathbf{H} \gets \mathbf{0}_{m \times m}$
    \State $\mathbf{e}_i \gets \mathbf{0}_{m}$ \Comment{temporary basis vector}
    \State $\mathbf{e}_j \gets \mathbf{0}_{m}$ \Comment{temporary basis vector}

    \Statex
    \Comment{--- Diagonal components ---}
    \For{$i \gets 1$ \textbf{to} $m$}
      \State $\mathbf{e}_i[i] \gets 1$
      \State $\mathbf{H}[i,i] \gets d_{2f\mathbf{q}}(\mathbf{e}_i)$
      \State $\mathbf{e}_i[i] \gets 0$
    \EndFor

    \Statex
    \Comment{--- Off-diagonal components ---}
    \For{$i \gets 1$ \textbf{to} $m$}
      \State $\mathbf{e}_i[i] \gets 1$
      \State $\mathbf{e}_j \gets \mathbf{0}_{m}$
      \For{$j \gets i+1$ \textbf{to} $m$}
        \State $\mathbf{e}_j[j] \gets 1$
        \State $h_{ij} \gets \frac{1}{2}\Big(
        d_{2f\mathbf{q}}(\mathbf{e}_i + \mathbf{e}_j)
        - \mathbf{H}[i,i] - \mathbf{H}[j,j]\Big)$
        \State $\mathbf{H}[i,j] \gets h_{ij}$;
        \quad $\mathbf{H}[j,i] \gets h_{ij}$
        \State $\mathbf{e}_j[j] \gets 0$
      \EndFor
      \State $\mathbf{e}_i[i] \gets 0$
    \EndFor

    \State \Return $\mathbf{H}$
  \end{algorithmic}
\end{breakablealgorithm}
\vspace*{0.2cm}

Higher-order derivatives can be computed following an analogous
procedure, using the function $d_{nf\mathbf{q}}$ or its vector-valued
counterpart. In addition, the function
\texttt{SPD(fsd, Aindx, q)}, included in the archived
implementation~\cite{penunuri_2026_zDNAODD}, is provided for the
computation of mixed partial derivatives of arbitrary order.  This
function returns the partial derivative of the scalar-valued function
\texttt{fsd}, specified by the index matrix \texttt{Aindx}, evaluated
at the point \texttt{q}. The first column of \texttt{Aindx} specifies
the variable indices with respect to which differentiation is performed,
while the second column specifies how many times the derivative is taken
with respect to each corresponding variable.

For example, the code in Listing~\ref{lst:d5fpd} computes
\begin{align}
  f_{xxyzz}(\mathbf{q})
  = \frac{\partial^5 f(\mathbf{q})}
  {\partial x^2\, \partial y\, \partial z^2},
\end{align}
for $f(x,y,z) = \sin(xyz) + \cos(xyz)$, evaluated at
$\mathbf{q} = [0.1 + i,\, 0.2 + i,\, 0.3 + i]$.
\vspace*{0.2cm}
\begin{lstlisting}[
style=fortran-style,
caption={Computation of $f_{xxyzz}$},
label={lst:d5fpd}]
program main
  use config_mod
  use dualzn_mod
  use dir_der_mod
  use test_functions_mod
  implicit none

  integer, allocatable, dimension(:,:) :: indx_mat
  complex(prec), parameter :: ii = (0.0_prec,1.0_prec)
  complex(prec), dimension(3) :: q
  complex(prec) :: dnfs

  q = [0.1_prec + ii, 0.2_prec + ii, 0.3_prec + ii]
  allocate(indx_mat(3,2))
  ! [[1,2],[2,1],[3,2]] <--> {{x,2},{y,1},{z,2}}
  indx_mat(:,1) = [1,2,3]
  indx_mat(:,2) = [2,1,2]

  write(*,"(A)") "--SPD--"
  dnfs = SPD(fstest, indx_mat, q)
  write(*,*) "d5f/dx2dydz2:", dnfs
  deallocate(indx_mat)
end program main
\end{lstlisting}

\subsection{Kinematic quantities: Position, velocity, acceleration,
  jerk, snap, and higher-order kinematic quantities}
\label{KinQsSec}

The computation of kinematic quantities is of central importance in the
fields of mechanisms and robotics~\cite{Gosselin1990,Merlet2006,ArikawaKeisuke}.
Among the classical approaches for this purpose, screw
theory~\cite{Rico1999,Custodio2017,Antonov2024} provides a powerful
framework for the computation of kinematic quantities. Its application,
however, typically requires an explicit description of the underlying
kinematic chain together with the identification of the corresponding
joint screws, and the extension to higher-order derivatives often leads
to expressions of considerable algebraic complexity. The methodology
proposed in this work instead provides a unified and systematic
procedure for computing kinematic quantities of arbitrary order directly
from the position vector $\mathbf{r}=\mathbf{f}(\mathbf{q})$, by
expressing them in terms of directional derivatives evaluated through
dual numbers.

Let $\mathbf{r} = \mathbf{f}(\mathbf{q})$ be a vector-valued position
function with $\mathbf{q} = \mathbf{q}(t)$, and let
$\mathbf{q}^{(j)} = d^j \mathbf{q}/dt^j$ denote its time derivatives.
The velocity, acceleration, jerk, and snap vectors can be
expressed in terms of vector directional derivatives as
\cite{PEONESCALANTE2024}:
\begin{align}
\mathbf{r}^{(1)} &= \mathbf{d}_1\!\left(\mathbf{q}^{(1)}\right), \label{veld1mlf}\\[4pt]
\mathbf{r}^{(2)} &= \mathbf{d}_2\!\left(\mathbf{q}^{(1)}\right) +
\mathbf{d}_1\!\left(\mathbf{q}^{(2)}\right), \label{aceld2mlf}\\[4pt]
\mathbf{r}^{(3)} &= \mathbf{d}_3\!\left(\mathbf{q}^{(1)}\right) + 
3\,\mathbf{d}_2\!\left(\mathbf{q}^{(1)},\mathbf{q}^{(2)}\right) +
\mathbf{d}_1\!\left(\mathbf{q}^{(3)}\right), \label{pulsod3mlf}\\[4pt]
\mathbf{r}^{(4)} &=  \mathbf{d}_4\!\left(\mathbf{q}^{(1)}\right) +
6\,\mathbf{d}_3\!\left((\mathbf{q}^{(1)})^{[2]},\mathbf{q}^{(2)}\right) + 
3\,\mathbf{d}_2\!\left(\mathbf{q}^{(2)}\right) +
4\,\mathbf{d}_2\!\left(\mathbf{q}^{(1)},\mathbf{q}^{(3)}\right) +
\mathbf{d}_1\!\left(\mathbf{q}^{(4)}\right), \label{jounced4mlf}
\end{align}
where, to lighten the notation, we define
\begin{align}
  \mathbf{d}_{k}&:=\mathbf{d}_{k\mathbf{f}\mathbf{q}}, \label{dk:dkfq}\\
  (\mathbf{q}^{(i)})^{[m]} &:=
  \underbrace{\mathbf{q}^{(i)},\ldots,\mathbf{q}^{(i)}}_{m},\label{qim:qiqim}\\
  \mathbf{d}_{k}(\mathbf{q}^{(i)}) &:= \mathbf{d}_{k}\!\left((\mathbf{q}^{(i)})^{[k]}\right).
  \label{dkqi:dkqik}
\end{align}

With this notation, the multivariate Fa\`a~di~Bruno expansion
\cite{BAEZA2017156} yields the $n$th total derivative of
$\mathbf{r}(\mathbf{q}(t))$ as
\begin{equation}
  \mathbf{r}^{(n)} =
  \sum_{\substack{m_1,\dots,m_n\ge0\\
  \sum j m_j = n}}
  \frac{n!}{\displaystyle
  \prod_{j=1}^{n} m_j!\,(j!)^{m_j}}
  \;
  \mathbf{d}_{m_1+\cdots+m_n}
  \!\left(
  \big(\mathbf{q}^{(1)}\big)^{[m_1]},
  \big(\mathbf{q}^{(2)}\big)^{[m_2]},\dots,
  \big(\mathbf{q}^{(n)}\big)^{[m_n]}
  \right).
  \label{eq:fdb-vector}
\end{equation}
This expression can be evaluated numerically, since we have an explicit
implementation for computing the required directional derivatives.
Nevertheless, Eq.~\eqref{eq:fdb-vector} can be reorganized by grouping
all terms that share the same number of directional arguments,
$k = m_1+\cdots+m_n$. This grouping parallels the structure of the
\emph{partial Bell polynomials} $B_{n,k}$ in the univariate
Fa\`a~di~Bruno formula, which collect all integer tuples
$(m_1,\dots,m_{n-k+1})$ satisfying
$\sum_{j=1}^{n-k+1} j\,m_j=n$
and
$\sum_{j=1}^{n-k+1} m_j=k$.
Although in the genuinely multivariate setting the derivatives
$\mathbf{d}_k$ act on vector arguments rather than scalar products, the
same combinatorial structure applies. Thus, for compactness, we
introduce the symbolic notation
\begin{equation}
  \mathbf{r}^{(n)} =
  \sum_{k=1}^{n}
  \mathbf{d}_k\!\left(
  \mathcal{B}_{n,k}\big(
  \mathbf{q}^{(1)},\mathbf{q}^{(2)},\dots,
  \mathbf{q}^{(n-k+1)}
  \big)
  \right),
  \label{eq:rnk-bell}
\end{equation}
where $\mathcal{B}_{n,k}$ denotes the symbolic Bell operator associated
with the scalar partial Bell polynomial $B_{n,k}$. Its role here is not
to define a scalar polynomial evaluated on vectors, but rather to encode
the admissible integer tuples $(m_j)$ together with their corresponding
Bell coefficients. In this interpretation, each monomial
\[
x_1^{m_1}\cdots x_{n-k+1}^{m_{n-k+1}}
\]
is mapped to the ordered list of repeated vector arguments
\[
(\mathbf{q}^{(1)})^{[m_1]},\dots,
(\mathbf{q}^{(n-k+1)})^{[m_{n-k+1}]},
\]
which are then supplied to the multilinear differential
$\mathbf{d}_k$.

The equivalence between
Eqs.~\eqref{eq:fdb-vector} and~\eqref{eq:rnk-bell}
follows by grouping the terms of Eq.~\eqref{eq:fdb-vector}
according to the total multiplicity
$k=\sum m_j$, which yields the classical definition of the
partial Bell polynomials,
\[
B_{n,k}(x_1,\ldots,x_{n-k+1})
=
\!\!\sum_{\substack{
m_1+\cdots+m_{n-k+1}=k\\
\sum j\,m_j=n}}
\frac{n!}
{\prod_{j=1}^{n-k+1} m_j!\,(j!)^{m_j}}
x_1^{m_1}\cdots x_{n-k+1}^{m_{n-k+1}}.
\]
The symbolic operator $\mathcal{B}_{n,k}$ therefore preserves the same
integer-partition structure as the scalar Bell polynomial, while the
actual evaluation is performed through the multilinear action of
$\mathbf{d}_k$ on vector arguments.

In practice, the coefficients and index pairs that define each term in
Eq.~\eqref{eq:rnk-bell} can be obtained in \textsc{Mathematica} using the
built-in \texttt{BellY} function, which symbolically generates the
partial Bell polynomials. The auxiliary routine \verb+FDBPairs+
(Listing~\ref{lst:math1}) automates this process by enumerating all
integer tuples $(m_j)$ together with the corresponding coefficients and
indices required to assemble $\mathbf{r}^{(n)}$ for arbitrary order.
Table~\ref{ejr4Bnk} illustrates the resulting terms for the case
$\mathbf{r}^{(4)}$.
\begin{lstlisting}[
    caption={Function to generate the coefficients and index
      pairs required to build $\mathbf{r}^{(n)}$.},
    label={lst:math1}]
(* Returns: {{coeff,k},{{j1,m1},{j2,m2},...}} with m_j>0 *)
FDBPairs[n_Integer?Positive] :=
  Module[{vars, cr},
   Reap[
     Do[
       vars = Array[x, n - k + 1];
       cr = CoefficientRules[BellY[n, k, vars], vars];
       Scan[Function[rule,
         With[{exps = First@rule, coeff = Last@rule, len = Length@vars},
           Sow@{{coeff, k},
             Pick[Transpose[{Range[len], exps}], exps, _?(# > 0 &)]}
         ]],
        cr],
       {k, 1, n}]
   ][[2, 1]]];

FDBPairs[4] // Reverse
(* {{{1, 4}, {{1, 4}}},
    {{6, 3}, {{1, 2}, {2, 1}}},
    {{3, 2}, {{2, 2}}},
    {{4, 2}, {{1, 1}, {3, 1}}},
    {{1, 1}, {{4, 1}}}} *)
\end{lstlisting}

\begin{table}[htb]
  \centering
  \caption{Terms of $\mathbf{r}^{(4)}(t)$ in nested-list notation.
    Each entry corresponds to $\{\{\text{coeff},k\},\{\{j_1,m_1\},\ldots\}\}$.}
  \label{ejr4Bnk}
  \vspace{0.5ex}
  \begin{tabular}{lll}
    \toprule
    \verb|{coeff,k}| & \verb|{{j,m},...}| & Directional-derivative term \\
    \midrule
    \verb|{1,4}| & \verb|{{1,4}}| &
      $\mathbf{d}_4 \big(\mathbf{q}^{(1)}\big)$ \\[4pt]
    \verb|{6,3}| & \verb|{{1,2},{2,1}}| &
      $6\,\mathbf{d}_3 \big((\mathbf{q}^{(1)})^{[2]},\mathbf{q}^{(2)}\big)$ \\[4pt]
    \verb|{3,2}| & \verb|{{2,2}}| &
      $3\,\mathbf{d}_2 \big(\mathbf{q}^{(2)}\big)$ \\[4pt]
    \verb|{4,2}| & \verb|{{1,1},{3,1}}| &
      $4\,\mathbf{d}_2 \big(\mathbf{q}^{(1)},\mathbf{q}^{(3)}\big)$ \\[4pt]
    \verb|{1,1}| & \verb|{{4,1}}| &
      $\mathbf{d}_1 \big(\mathbf{q}^{(4)}\big)$ \\
    \bottomrule
  \end{tabular}
\end{table}

\paragraph{Computational correspondence}
Each element of \texttt{FDBPairs(n)} represents one term in
Eq.~\eqref{eq:fdb-vector}, where \texttt{coeff} corresponds to the
combinatorial factor
$n! / \prod m_j!\,(j!)^{m_j}$,
and \texttt{\{j,m\}} encodes the multiplicities $m_j$.
Defining the map
\[
\Phi\!\left(\big\{\{\texttt{coeff},k\},\{\{j_1,m_{j_1}\},\dots\}\big\}\right)
:= \texttt{coeff}\,
\mathbf{d}_{k}\!\Big(
(\mathbf{q}^{(j_1)})^{[m_{j_1}]},\dots\Big),
\]
we obtain
\begin{equation}
\mathbf{r}^{(n)} = \sum_{\tau\in \texttt{FDBPairs}(n)} \Phi(\tau),
\label{eq:sum-FDBPairs}
\end{equation}
which reproduces the same combinatorial and multilinear expansion as
Eqs.~\eqref{eq:fdb-vector} and~\eqref{eq:rnk-bell}.

Although the symbolic construction based on Eq.~\eqref{eq:rnk-bell} is
convenient for analytical purposes, it is not strictly necessary in a
computational setting. The same combinatorial structure can be generated
programmatically by enumerating all admissible integer tuples
$(m_1,\dots,m_n)$ satisfying $\sum j\,m_j = n$, and then evaluating the
corresponding directional-derivative terms numerically. In other words,
the symbolic machinery of $B_{n,k}$ (\texttt{BellY}) can be bypassed
altogether: once the coefficients and index pairs $\{\{j,m_j\}\}$ have
been generated (as in the output of \texttt{FDBPairs}), the quantities
$\mathbf{r}^{(n)}$ can be assembled using standard arrays and iterative
loops (or recursion, for convenience). This approach is particularly
suitable for low-level languages such as Fortran, where numerical
efficiency and memory locality are essential. The Fortran module
\verb+fdb_bell_mod.f90+ in \cite{penunuri_2026_zDNAODD} includes the function
\verb+generate_fdb_bell+, which generates the coefficients and index
pairs required to build $\mathbf{r}^{(n)}$ without resorting to any
symbolic computation. Listing~\ref{lst:fortraFDBC} presents an example of
its use.

\vspace*{0.2cm}
\begin{lstlisting}[
    style=fortran-style,
    caption={Fortran program to compute coefficients and indices for
  $\mathbf{r}^{(n)}(t)$.},
label={lst:fortraFDBC}]
program main
  use fdb_bell_mod
  implicit none
  integer :: n
  type(term_t), allocatable :: terms(:)

  write(*,'(a)') 'nd:'
  read(*,*) n

  terms = generate_fdb_bell(n)
  call print_terms(n, terms)
  deallocate(terms)
end program main


!----- output for nd = 5 -----
!nd:
!5
!--- Faa di Bruno Bell Terms for n = 5 ---
!{[1,1],[[5,1]]}
!{[5,2],[[1,1],[4,1]]}
!{[10,2],[[2,1],[3,1]]}
!{[10,3],[[1,2],[3,1]]}
!{[15,3],[[1,1],[2,2]]}
!{[10,4],[[1,3],[2,1]]}
!{[1,5],[[1,5]]}
\end{lstlisting}

\begin{figure}[htb]
  \centering
  \includegraphics[scale=0.3]{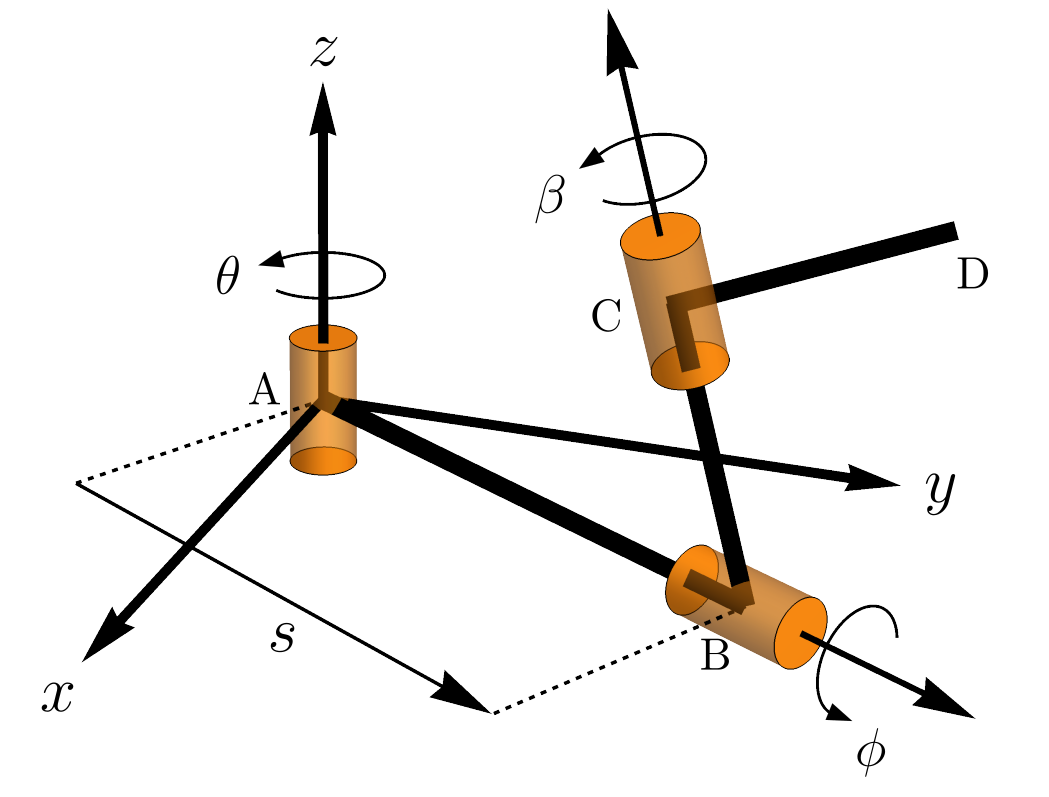}
  \caption{RCR robot manipulator benchmark example.}
  \label{RCRFig}
\end{figure}

The RCR robot manipulator shown in Fig.~\ref{RCRFig} is reconsidered as
a benchmark example to illustrate the applicability of the proposed
formulation to arbitrary-order kinematic quantities. The same geometric
parameters and generalized-coordinate values previously reported
in~\cite{PEONESCALANTE2024} are adopted here, allowing direct comparison
with the lower-order results presented therein. Although the proposed
framework applies to derivatives of arbitrary order, only the
fifth-order kinematic quantity is reported explicitly in the present
example.

Table~\ref{KQJDJ} reports the resulting fifth-order kinematic quantity
for the end effector of the manipulator corresponding to
$\mathbf{q}^{(5)} = [1~~3~~5~~7]^{\text{T}}$.

\begin{table}[htb]
  \caption{Fifth-order kinematic quantity for the end effector of the
    RCR robot manipulator.}
\centering
\begin{tabular}{lr}
\toprule
 Kinematic quantity & Component values  \\
\midrule
                      & $5830.00$\\
$\mathbf{r}_D^{(5)}$ ~\normalsize{[L]/[t]$^5$} & $-1020.00$\\
\multirow{-3}{*}{}                            & $-2735.00$\\
\bottomrule
\end{tabular}\label{KQJDJ}
\end{table}

\subsection{Taylor series method for solving ODEs}
The Taylor Series Method (TSM)~\cite{CHANG1994209,Nedialkov2005,BARRIO2011}
provides a conceptually simple and mathematically rigorous framework for
approximating the solutions of systems of ordinary differential equations
(ODEs). In practice, however, its implementation is often hindered by the
need to compute higher-order derivatives of the right-hand side function.
As a result, most existing studies have been limited either to particular
cases or to symbolic differentiation performed using computer algebra
systems such as \textsc{Mathematica}. With the availability of a robust,
efficient, and transparent procedure for evaluating \emph{time}
derivatives, the application of the TSM becomes straightforward.
Although this approach may not match the computational efficiency of
simpler algorithms such as the classical fourth-order Runge--Kutta
method, its ease of generalization and intrinsic accuracy make it a
competitive and reliable option for the numerical integration of ODE
systems.

Let the system of ordinary differential equations be defined as
\begin{equation}
  \label{eq:ode_system}
  \begin{aligned}
    \mathbf{x}^{(1)}(t) &= \mathbf{F}(\mathbf{q}(t)), \\
    \mathbf{x}(t_0) &= \mathbf{x}_0,
  \end{aligned}
\end{equation}
where
\begin{equation}
  \begin{aligned}
    \mathbf{q}(t) &= [\,x_1(t),\, x_2(t),\, \dots,\, x_n(t),\, t\,], \\
    \mathbf{F}(\mathbf{q}) &= [\,f_1(\mathbf{q}),\, f_2(\mathbf{q}),\,
    \dots,\, f_n(\mathbf{q})\,],
  \end{aligned}
\end{equation}
with each function $f_i:\mathbb{R}^{n+1} \to \mathbb{R}$ representing the
$i$th component of the vector field. Here, $\mathbf{x}_0$ denotes the
vector of initial conditions, $t$ is the independent variable, and, as
before, $\mathbf{x}^{(k)}$ denotes the $k$th derivative with respect to $t$.

According to a fourth-order TSM, the solution at $t_0 + h$ is given by
\begin{align}
  \mathbf{x}(t_0 + h) = \mathbf{x}(t_0)
  + h\,\mathbf{x}^{(1)}(t_0)
  + \frac{h^2}{2}\,\mathbf{x}^{(2)}(t_0)
  + \frac{h^3}{3!}\,\mathbf{x}^{(3)}(t_0)
  + \frac{h^4}{4!}\,\mathbf{x}^{(4)}(t_0).
\end{align}
With the proposed formulation, all required derivatives are
straightforward to compute. A pseudocode description of the method is
presented below.

\begin{breakablealgorithm}
  \caption{Pseudocode for \textsc{TSMDD}$(\mathbf{F},\,\mathbf{x}_0,\,\mathbf{t})$}
  \label{alg:tsmdd}
  \begin{algorithmic}[1]
    \Require Procedure $\mathbf{F}:\mathbb{D}^{n+1} \to \mathbb{D}^{n}$
    \Require $\mathbf{x}_0 \in \mathbb{R}^{n}$ \Comment{initial conditions}
    \Require $\mathbf{t} = [t_1,\dots,t_m] \in \mathbb{R}^{m}$ \Comment{time vector}
    \Ensure $\mathbf{X} \in \mathbb{R}^{m \times n}$ \Comment{solution matrix (excluding $t$)}

    \State $n \gets \text{size}(\mathbf{x}_0)$;
    \quad $m \gets \text{size}(\mathbf{t})$;
    \quad $n_{\text{int}} \gets m - 1$
    \State $\mathbf{X}[1,:] \gets \mathbf{x}_0$
    \State $\mathbf{h} \gets [\,t_2 - t_1,\, t_3 - t_2,\, \dots,\, t_m - t_{m-1}\,]$

    \State Initialize $\mathbf{q}^{(1)} = [\,\mathbf{0}_{1\times n},\,1\,]$
    \State Initialize $\mathbf{q}^{(2)} = \mathbf{0}_{1\times (n+1)}$
    \State Initialize $\mathbf{q}^{(3)} = \mathbf{0}_{1\times (n+1)}$

    \For{$i \gets 1$ \textbf{to} $n_{\text{int}}$}
      \State $t_i \gets \mathbf{t}[i]$;\quad $h \gets \mathbf{h}[i]$
      \State $\mathbf{q} \gets [\,\mathbf{X}[i,:],\,t_i\,]$

      \State $\mathbf{x}^{(1)} \gets \mathbf{F}(\mathbf{q})$
      \State $\mathbf{q}^{(1)}[1{:}n] \gets \mathbf{x}^{(1)}$

      \State $\mathbf{x}^{(2)} \gets
      \mathbf{d}_{1\mathbf{Fq}}(\mathbf{q}^{(1)})$
      \State $\mathbf{q}^{(2)}[1{:}n] \gets \mathbf{x}^{(2)}$

      \State $\mathbf{x}^{(3)} \gets
      \mathbf{d}_{2\mathbf{Fq}}(\mathbf{q}^{(1)})
      + \mathbf{d}_{1\mathbf{Fq}}(\mathbf{q}^{(2)})$
      \State $\mathbf{q}^{(3)}[1{:}n] \gets \mathbf{x}^{(3)}$

      \State $\mathbf{x}^{(4)} \gets
      \mathbf{d}_{3\mathbf{Fq}}(\mathbf{q}^{(1)})
      + 3\,\mathbf{d}_{2\mathbf{Fq}}(\mathbf{q}^{(1)},\mathbf{q}^{(2)})
      + \mathbf{d}_{1\mathbf{Fq}}(\mathbf{q}^{(3)})$

      \State $\mathbf{X}[i{+}1,:] \gets \mathbf{X}[i,:]
      + h\,\mathbf{x}^{(1)}
      + \dfrac{h^2}{2}\,\mathbf{x}^{(2)}
      + \dfrac{h^3}{6}\,\mathbf{x}^{(3)}
      + \dfrac{h^4}{24}\,\mathbf{x}^{(4)}$
    \EndFor

    \State \Return $\mathbf{X}$
  \end{algorithmic}
\end{breakablealgorithm}

\begin{figure}[htb]
  \centering
  \includegraphics[scale=0.6]{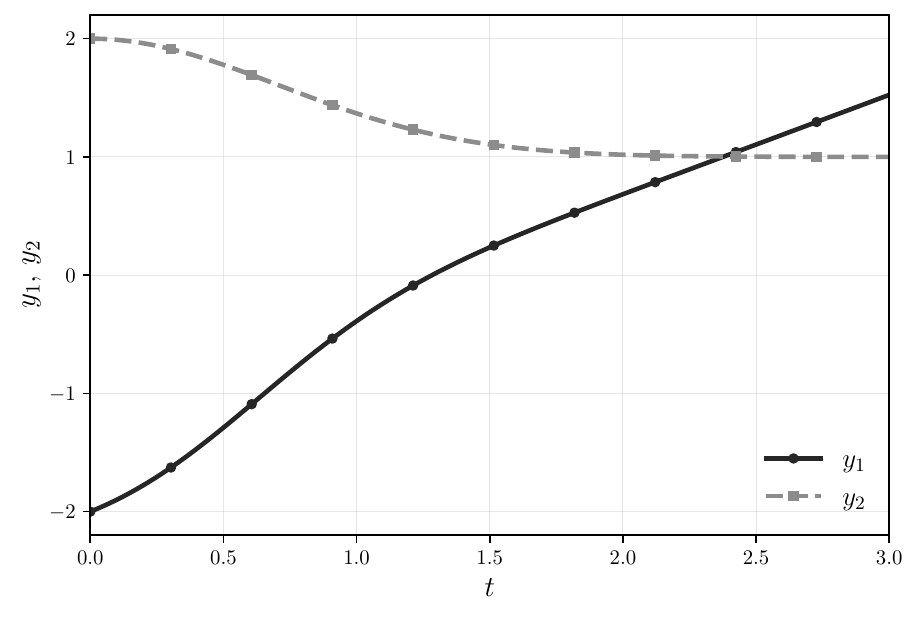}
  \caption{Solution obtained using the Taylor Series Method,
  Algorithm~\ref{alg:tsmdd}.}
  \label{TSMFig}
\end{figure}

Figure~\ref{TSMFig} shows the numerical solution of the following system
of ordinary differential equations:
\begin{align}
  y_1'(t) &= 2\,t\,\exp(-t^2)\bigl(1 - t\cos y_2 \bigr) + \sin y_2, \\
  y_2'(t) &= -2\,t\,\exp(-t^2),
\end{align}
subject to the initial conditions
\begin{align}
  y_1(0) &= -2, \\
  y_2(0) &= 2.
\end{align}
The system is solved using a Fortran implementation of
Algorithm~\ref{alg:tsmdd}.

\section{Conclusions}
This paper has introduced a dual-number-based formulation for
computing directional derivatives of arbitrary order along multiple
directions.  The main result is a practical mechanism for evaluating
the general symmetric multilinear form
$\mathbf{d}_n(\mathbf{v}_1,\ldots,\mathbf{v}_n)$ using only
evaluations of $\mathbf{d}_n(\mathbf{v})$, obtained from a single
function call with appropriately constructed dual variables. From the
coordinate representation of these multilinear forms, mixed
directional and mixed partial derivatives can be recovered without
symbolic differentiation and without assembling full higher-order
derivative tensors.

From a computational standpoint, the proposed approach is particularly
advantageous in regimes where the explicit enumeration of partial
derivatives becomes combinatorially prohibitive, but where higher-order
derivatives are required only along selected directions—a common
situation in sensitivity analysis and in problems involving higher-order
effects. The method was demonstrated on a high-dimensional benchmark
function, showing that higher-order directional derivatives can be
computed efficiently even for problems involving thousands of
variables. The proposed approach is not intended as a replacement for
full tensor-based higher-order automatic differentiation, but rather as
a targeted tool for efficiently computing selected higher-order
directional derivatives along prescribed directions. In addition, the
availability of this functionality in \textsc{Fortran}, where comparable
higher-order directional differentiation tools are currently limited,
makes the proposed implementation particularly relevant for large-scale
scientific computing and number-crunching applications.

A second contribution of this work is the application of the framework
to the systematic computation of arbitrary-order kinematic quantities
from a position mapping $\mathbf{r}=\mathbf{f}(\mathbf{q})$. By
expressing total time derivatives in terms of vector directional
derivatives and multivariate Fa\`a~di~Bruno combinatorics—organized
through Bell-polynomial index structures—the method provides a unified
procedure for computing velocity, acceleration, jerk, snap, and
higher-order kinematic quantities. The RCR manipulator example
illustrates the extension of the framework to higher-order
kinematic derivatives.

Finally, the availability of reliable high-order time derivatives
enables a transparent implementation of the Taylor Series Method for the
numerical solution of systems of ordinary differential equations. While
the resulting integrator is not necessarily the most cost-effective
option when compared with classical low-order schemes, it offers a
rigorous and easily generalizable alternative when high-order expansions
are desirable. The complete \textsc{Fortran} implementation accompanying
this work is publicly available, facilitating reuse, verification, and
further extension of the proposed approach.





\section*{Declaration of competing interest}
The authors declare that they have no known competing financial 
interests or personal relationships that could have appeared to
influence the work reported in this paper.

\section*{Data availability statement}
The source code supporting the findings of this study is publicly
available in a GitHub repository (\url{https://github.com/fpenunuri/DNAODD}) and has
been archived with a persistent identifier at
Zenodo~\cite{penunuri_2026_zDNAODD}. The archived version corresponds
to the implementation used to generate the results reported in this
paper.





\end{document}